\documentclass[journal]{IEEEtran}
\usepackage{mathrsfs}
\usepackage{graphicx}
\usepackage{latexsym}
\usepackage{algorithm}
\usepackage{algorithmic}
\usepackage{subfigure}
\usepackage{amsmath}
\usepackage{cite}
\usepackage{amssymb}
\usepackage{multirow}
\usepackage{amsthm}
\usepackage{color}
\usepackage{diagbox}
\usepackage{wrapfig}

\begin{document}
\title{Convergence Analysis of a Cooperative Diffusion Gauss-Newton Strategy}

\author{Mou Wu, Naixue Xiong, Liansheng Tan
\thanks{Mou Wu is with College of Intelligence and Computing, Tianjin University, Tianjin, 300350, PR China, and also with School of Computer Science and Technology, Hubei University of Science and Technology, Xianning 437100, PR China. (Email: mou.wu@163.com).}

\thanks{Naixue Xiong is with College of Intelligence and Computing, Tianjin University, Tianjin, 300350, PR China. (Email: xiongnaixue@gmail.com)}

\thanks{Liansheng Tan is with Dept. of Computer Science, Central China Normal University, Wuhan 430079, PR China. (Email: l.tan@mail.ccnu.edu.cn)}
}

\maketitle
\begin{abstract}
In this paper, we investigate the convergence performance of a cooperative diffusion Gauss-Newton (GN) method, which is widely used to solve the nonlinear least squares problems (NLLS) due to the low computation cost compared with Newton's method. This diffusion GN collects the diversity of temporal-spatial information over the network, which is used on local updates.  In order to address the challenges on convergence analysis, we firstly consider to form a global recursion relation over spatial and temporal scales since the traditional GN is a time iterative method and the network-wide NLLS need to be solved. Secondly, the derived recursion related to the network-wide deviation between the successive two iterations is ambiguous due to the uncertainty of  descent discrepancy in GN update step between two versions of cooperation and non-cooperation. Thus, an important work is to derive the boundedness conditions of this discrepancy. Finally, based on the temporal-spatial recursion relation and the steady-state equilibria theory for discrete dynamical systems, we obtain the sufficient conditions for algorithm convergence, which require the good initial guesses, reasonable step size values and network connectivity. Such analysis provides a guideline for the applications based on this diffusion GN method.

\end{abstract}


\section{Introduction}
\IEEEPARstart{G}{auss}-Newton method has found wide applications, such as deep learning in artificial intelligence and neural network \cite{Schmidhuber2015Deep,wang2018deep}, and parameter estimate in a networked system \cite{Zhao2007Information,Xiao2013Convergence,Kaur2017A}. Deriving from Newton's method, GN algorithm discards the second-order terms in the computation of Hessian for small residual NLLS problems, thereby resulting in saving in computation. Such the amount of computations can be further reduced via the mathematical process.  In order to compute easily the first derivative of objective function, the perturbed GN method is proposed in \cite{Gratton2007Approximate}, where a perturbed derivative version substitutes the original one. The truncated GN method \cite{Bao2017Approximate} is proposed to implement the inexact update instead of exact one. The truncated-perturbed GN method \cite{Bao2017Approximate} integrates the above two advantages into the update step.

Many scenarios can be modeled as the NLLS problem depended on the performance of GN, such as computer vision \cite{Xiong2014Supervised}, image alignment and reconstruction \cite{Torre2008Parameterized,Schweiger2005Gauss}, network-based localization \cite{souza2016target,bejar2011distributed}, signal processing for direction-of-arrival estimation and frequency estimation  \cite{Jensen2013Nonlinear}, logistic regression \cite{Zhuang2015distributed} and power system state estimation \cite{Xiao2013Convergence,cosovic2016distributed,Minot2016A}. Despite the widespread utility, it is difficult for exploiting the original GN method as a fully cooperative scheme for a distributed network, since its iteration rule involves the matrix inverse operator, which is ideally suited to be implemented in a centralized way. However, for the well known advantages such as load balancing and robustness, distributed algorithm with the improvement of performance is preferred.

The purpose of this work is to analyze the convergence of  a cooperative diffusion GN strategy over a distributed network, where every node sense the temporal data that is variable over the spatial domain. Several diffusion GN methods \cite{Calafiore2010A,Bejar2015Distributed} are proposed for solving the localization problem in wireless sensor networks. However, they are centralized in nature and implemented in a non-cooperative way, in which the local intermediate estimates are not shared over the diffusion network.

\textbf{Notation:}  The operator $(\cdot)^{T}$ denotes the transpose for matrix or vector, the operator $(\cdot)^{-1}$ denotes the inverse of a non-singular matrix. The capital letters are used when the matrices are denoted, while the small letters are used when the vectors or scalars  are denoted. The Euclidean norm of a vector $x$ is written as $\|x\|$, 2-norm and \textit{Frobenius} norm of a matrix $G$ is denoted by $\|G\|$ and $\|G\|_{F}$, respectively. $I_{N}$ and $\textmd{1}_{N}$ denote the $N\times N$ identity matrix and $N\times 1$  vector whose every entry is 1, respectively. We will use subscripts $k$, $l$, $u$ and $t$ to denote node, and superscript $j$, $i$ to denote time.

\section{Description of cooperative diffusion Gauss-Newton solution}
\subsection{Centralized solution}

For an adaptive network represented by a set $\mathcal{N}=\{1,\cdots,N\}$ , we would like to estimate a $M\times1$ unknown parameter vector $x=[x_{1},\cdots,x_{M}]^{T}$ belonging to a closed convex set $\mathbb{X}$. Let $f(x)=[f_{1}(x),\cdots,f_{N}(x)]^{T}:\Re^{M}\longrightarrow\Re^{N}$ be a  continuous and differentiable global cost function throughout the network, where $f_{k}(x):\Re^{M}\longrightarrow\Re$ is the individual cost function associated with node $k\in \mathcal{N}$ by collecting the measurements from the related events. The estimation problem can be formulated as

\begin{equation}  \label{eq1}
\mathop{min}\limits_{x}\parallel f(x)\parallel^{2}.
\end{equation}

By rewriting $\|f(x)\|^{2}=\sum^{N}_{k=1}|f_{k}(x)|^{2}$, the object of each node in the network is to seek a $M\times1$  vector $x$ that solve the following Non-Linear Least Squares (NLLS) problem with the form

\begin{equation}\label{eq2}
 \mathop{min}\limits_{x}\sum\limits_{k=1}^{N}|f_{k}(x)|^{2}.
\end{equation}

The GN method is well recognized for solving NLLS problems. Let us consider a fusion center (FC) that can communicate with all nodes in the network. Given an initial good guess $x^{0}$, a centralized scheme can be implemented on FC based on the GN update rule in an iterative way

\begin{equation}\label{eq3}
 x^{i+1}=x^{i}-\alpha^{i}d^{i},
\end{equation}
where $x^{i}$ is the estimation of $x$ at iteration $i$, $d^{i}$ denotes a descent direction of GN, and $\alpha^{i}$ is the step size parameter that ensure $x^{i+1}$ is nearer a stationary point than $x^{i}$.

In this paper, we adopt the following assumption for the above optimization problem.

\textbf{Assumption 1.}

(1) The stationary points $x^{s}\in\Re^{M}$ that satisfy
$$\nabla f(x^{s})=2F^{T}(x^{s})f(x^{s})=0$$
always exist, where $F(x)$ is the Jacobian of $f(x)$ with the size $N\times M$ and the entries $F(x)_{k,l}=\partial f_{k}(x)/\partial x_{l}$, $1\leq k \leq N, 1\leq l \leq M$.

(2) The notations $\lambda _{min}(\cdot)$ and $\lambda _{max}(\cdot)$ are denoted as the minimum and maximum eigenvalues. For all $x\in \mathbb{X}$ and $k\in N$, let
$$\Sigma_{min}=\textrm{min}\sqrt{\lambda _{min}(F^{T}(x)F(x))}$$
$$\Sigma_{max}=\textrm{max}\sqrt{\lambda _{max}(F^{T}(x)F(x))},$$
where $0<\Sigma_{min}<\Sigma_{max}<\infty$.

Under Assumption 1, the approximate Hessian $F^{T}(x)F(x)$ of $f(x)$ is positive definite. Thereby, a local minimizer of $f(x)$ denoted by $x^{*}$ that belongs to the set of stationary points always exist \cite{kelley1999iterative,Bj1996Numerical}. Thus, the descent direction of GN update is written as

\begin{equation}\label{eq4}
d^{i}=[F^{T}(x^{i})F(x^{i})]^{-1}F^{T}(x^{i})f(x^{i}).
\end{equation}

By rewriting
\begin{equation}\label{eq5}
F(x)= \textrm{col}\{\frac{\partial f_{1}(x)}{\partial x},\frac{\partial f_{2}(x)}{\partial x},\cdots,\frac{\partial f_{N}(x)}{\partial x}\} \;\;\;(N\times M)
\end{equation}
and defining
\begin{equation}\label{eq6}
F_{k}(x)\triangleq \frac{\partial f_{k}(x)}{\partial x}, \;\;\;(1\times M)
\end{equation}
we get
\begin{equation}\label{eq7}
d^{i}=[\sum\limits_{k=1}^{N}F^{T}_{k}(x^{i})F_{k}(x^{i})]^{-1}\sum\limits_{k=1}^{N}F^{T}_{k}(x^{i})f_{k}(x^{i}).
\end{equation}

Therefore, we have the following GN iteration update
\begin{equation}\label{eq8}
 x^{i+1}=x^{i}-\alpha^{i}[\sum\limits_{k=1}^{N}F^{T}_{k}(x^{i})F_{k}(x^{i})]^{-1}\sum\limits_{k=1}^{N}F^{T}_{k}(x^{i})f_{k}(x^{i}).
\end{equation}

To successfully implement (\ref{eq8}) in a centralized way, we assume that the FC can communicate with all nodes over network and the same initial estimate is given by $x_{k}^{0}=x^{0}, k\in \mathcal{N}$. In the centralized GN algorithm,  the computation results of  $F^{T}_{k}(x^{i})F_{k}(x^{i})$ and $F^{T}_{k}(x^{i})f_{k}(x^{i})$ from each node $k$ are aggregated by the FC to obtain the new estimate $x^{i+1}$ based on (\ref{eq8}). Then the estimate $x^{i+1}$ is returned to all nodes until an appropriate termination condition is satisfied, for example $\|x^{i+1}-x^{i}\|\leq \varepsilon$ or $i=I$, where $\varepsilon$ and $I$ are the predefined minimum norm decline and the maximum number of iterations, respectively. Thus, the centralized GN includes actually a step of diffusion for new estimate $x^{i+1}$ form FC to individual nodes.

In this paper, we adopt the constant step size $\alpha_{k}^{i}=\alpha\in(0,1]$ for the subsequent development and analysis.

\subsection{Diffusion Gauss-Newton}

Consider the adaptive network $\mathcal{N}$, where any node $k$ at time $i$ receives a set of estimates $\{x_{l}^{i}\}_{l\in\mathcal{N}_{k}}$ from all its 1-hop neighbors $\mathcal{N}_{k}$ including itself. Thus,  the local estimates $\{x_{l}^{i}\}_{l\in\mathcal{N}_{k}}$ is combined in a weighted combination way denoted by
\begin{equation}\label{eq9}
 \mathcal{X}_{k}^{i}=\sum\limits_{l\in\mathcal{N}_{k}} c_{kl}x_{l}^{i},
\end{equation}
where $c_{kl}$ is the weighted coefficient between node $k$ and $l\in\mathcal{N}_{k}$. And the conditions

\begin{equation}\label{eq10}
\sum\limits_{l\in\mathcal{N}_{k}} c_{kl}=1\;\; and\;\; c_{kl}\in[0,1]\;\; for\;\; l\in\mathcal{N}_{k}
\end{equation}
is satisfied.

Once the aggregate estimate $\mathcal{X}_{k}^{i}$ is obtained as the local weighted estimate, any node $k$ in the network can implement the GN update step as follows:

\begin{equation}\label{eq11}
 x_{k}^{i+1}=\mathcal{X}_{k}^{i}-\alpha[Q_{k}^{i}(\mathcal{X})]^{-1}q_{k}^{i}(\mathcal{X}),
\end{equation}
where we define
\begin{equation}\label{eq12}
\begin{aligned}
 Q_{k}^{i}(\mathcal{X})&\triangleq F^{T}_{l\in\mathcal{N}_{k}}(\mathcal{X}_{l}^{i})F_{l\in\mathcal{N}_{k}}(\mathcal{X}_{l}^{i})\\ &\triangleq\sum\limits_{l\in\mathcal{N}_{k}}F^{T}_{l}(\mathcal{X}_{l}^{i})F_{l}(\mathcal{X}_{l}^{i})
\end{aligned}
\end{equation}
and
\begin{equation}\label{eq13}
\begin{aligned}
 q_{k}^{i}(\mathcal{X})&\triangleq F^{T}_{l\in\mathcal{N}_{k}}(\mathcal{X}_{l}^{i})f_{l\in\mathcal{N}_{k}}(\mathcal{X}_{l}^{i})\\
 &\triangleq \sum\limits_{l\in\mathcal{N}_{k}}F^{T}_{l}(\mathcal{X}_{l}^{i})f_{l}(\mathcal{X}_{l}^{i}).
\end{aligned}
\end{equation}

Removing the aggregate step of diffusion GN algorithm, we obtain a non-cooperative diffusion GN algorithm, where each node in the network acts as the FC to implement the centralized GN by communicating with all immediate neighbors. Its GN update step is given by

\begin{equation}\label{eq111}
 x_{k}^{i+1}=x_{k}^{i}-\alpha[Q_{k}^{i}(x_{k}^{i})]^{-1}q_{k}^{i}(x_{k}^{i}),
\end{equation}
where we define
\begin{equation}\label{eq121}
\begin{aligned}
 Q_{k}^{i}(x_{k}^{i})&\triangleq F^{T}_{l\in\mathcal{N}_{k}}(x_{k}^{i})F_{l\in\mathcal{N}_{k}}(x_{k}^{i})\\ &\triangleq\sum\limits_{l\in\mathcal{N}_{k}}F^{T}_{l}(x_{k}^{i})F_{l}(x_{k}^{i})
\end{aligned}
\end{equation}
and
\begin{equation}\label{eq131}
\begin{aligned}
 q_{k}^{i}(x_{k}^{i})&\triangleq F^{T}_{l\in\mathcal{N}_{k}}(x_{k}^{i})f_{l\in\mathcal{N}_{k}}(x_{k}^{i})\\
 &\triangleq \sum\limits_{l\in\mathcal{N}_{k}}F^{T}_{l}(x_{k}^{i})f_{l}(x_{k}^{i}).
\end{aligned}
\end{equation}
Note that the expression on arguments in (\ref{eq12}) (\ref{eq13}) (\ref{eq121}) (\ref{eq131}) shows the main difference between cooperative and non-cooperative algorithms.

The question that remains is how well does the diffusion GN algorithm perform in terms of its expected convergence behavior. First, what are the sufficient conditions of convergence for the diffusion GN algorithm? Second, is better the diffusion GN algorithm on convergence, compared with its non-cooperative counterpart? In other words, what are the benefits of cooperation? The following analysis and simulations will answer the above questions.

\section{Convergence analysis}

\subsection{Assumptions and data model}

To proceed the analysis, several reasonable assumptions need to be given as is commonly done in the literature \cite{Xiao2013Convergence,Johnson1991Topics}.

\textbf{Assumption 2.}

(1) $f_{l\in\mathcal{N}_{k}}(x_{k}^{i})$ is bounded for all $x_{k}^{i}\in\mathbb{X}\subset \mathbb{R}^{M}$ near  $x^{*}$, and satisfies
$$\|f_{l\in\mathcal{N}_{k}}(x_{k}^{i})\|\leq e_{max}$$
and
$$\|f_{l\in\mathcal{N}_{k}}(x^{*})\|= e_{min},$$
where $\|f_{l\in\mathcal{N}_{k}}(x^{*})\|$ denotes the  minimum value of $\|f_{l\in\mathcal{N}_{k}}(x_{k}^{i})\|$ when evaluated at $x_{k}^{i}=x^{*}$.

(2) For all $x\in \mathbb{X}$ and $k=1,\ldots,N$, let

$$\sigma_{min}=\textrm{min}\sqrt{\lambda _{min}(F_{k}^{T}(x)F_{k}(x))}$$
and
$$\sigma_{max}=\textrm{max}\sqrt{\lambda _{max}(F_{k}^{T}(x)F_{k}(x))},$$
where $0<\sigma_{min}<\sigma_{max}<\infty.$

(3) Both $F_{l\in\mathcal{N}_{k}}(x)$ and $F_{k}(x)$ are Lipschitz continuous on $\mathbb{X}$ with Lipschitz constant $\omega>0$ such that
$$\|F_{l\in\mathcal{N}_{k}}(x)-F_{l\in\mathcal{N}_{k}}(y)\|\leq \omega\|x-y\|$$
and
$$\|F_{k}(x)-F_{k}(y)\|\leq \omega\|x-y\|$$
for all $x,y\in \mathbb{X}$. Furthermore, we have the following results \cite{Eriksson2004Applied}
$$\|F_{k}^{T}(x)f_{k}(x)-F_{k}^{T}(y)f_{k}(y)\| \leq \gamma_{f}\|x-y\|$$
and
$$\|F_{k}^{T}(x)F_{k}(x)-F_{k}^{T}(y)F_{k}(y)\| \leq \gamma_{F}\|x-y\|,$$
where $\gamma_{f}\geq\omega(e_{max}+\Sigma_{max})$ and $\gamma_{F}\geq2\Sigma_{max}\omega$ are the corresponding Lipschitz constants.

In addition, the studying of the local convergence behavior need to be considered from the global view of network, since the performance of individual node depends on  the whole network including  cooperation rule and  network topology. Thus, we introduce the global quantities

$$x_{G}^{i}\triangleq \textrm{col}\{x_{1}^{i},\ldots,x_{N}^{i}\},\;\;\;(N M\times1)$$
$$\mathcal{X}_{G}^{i}\triangleq \textrm{col}\{\mathcal{X}_{1}^{i},\ldots,\mathcal{X}_{N}^{i}\},\;\;\;(N M\times1)$$
$$\overline{x}^{*}\triangleq \textrm{col}\{x^{*},\ldots,x^{*}\},\;\;\;(N M\times1)$$
$$D_{G}^{i}\triangleq \textrm{col}\{D_{1}^{i},\ldots,D_{N}^{i}\},\;\;\;(N M\times1)$$
$$d_{G}^{i}\triangleq \textrm{col}\{d_{1}^{i},\ldots,d_{N}^{i}\},\;\;\;(N M\times1)$$
where
$$D_{k}^{i}\triangleq [Q_{k}^{i}(\mathcal{X})]^{-1}q_{k}^{i}(\mathcal{X}),\;\;\;k\in\mathcal{N},$$
and
$$d_{k}^{i}\triangleq [Q_{k}^{i}(x_{k}^{i})]^{-1}q_{k}^{i}(x_{k}^{i}),\;\;\;k\in\mathcal{N},$$
$$A(x^{i}_{G})\triangleq \textrm{diag}\{F_{l\in\mathcal{N}_{1}}(x_{1}^{i}),\ldots,F_{l\in\mathcal{N}_{N}}(x_{N}^{i})\},\;\;\;(NN\times NM)$$
$$A(\overline{x}^{*})\triangleq \textrm{diag}\{F_{l\in\mathcal{N}_{1}}(x^{*}),\ldots,F_{l\in\mathcal{N}_{N}}(x^{*})\},\;\;\;(NN\times NM)$$
$$b(x^{i}_{G})\triangleq \textrm{col}\{f_{l\in\mathcal{N}_{1}}(x_{1}^{i}),\ldots,f_{l\in\mathcal{N}_{N}}(x_{N}^{i})\},\;\;\;(NN\times 1)$$
$$b(\overline{x}^{*})\triangleq \textrm{col}\{f_{l\in\mathcal{N}_{1}}(x^{*}),\ldots,f_{l\in\mathcal{N}_{N}}(x^{*})\},\;\;\;(NN\times 1)$$
where $\textrm{diag}(\cdot)$  is a block diagonal matrix whose entries are those of the column vector $\{\cdot\}$.

An $N\times N$ aggregate matrix $C$ can be given with non-negative real entries $\{c_{kl}\}$ that is redefined with the following conditions

\begin{equation}\label{eq14}
 c_{kl}=0\;\; \textrm{if}\;\;l\notin\mathcal{N}_{k}\;\;  \textrm{and}\;\;  \sum\limits_{l=1}^{N}c_{kl}=1,c_{kl}\geq0.
\end{equation}

Conditions (\ref{eq14}) indicate that the sum of all entries on each row of the matrix $C$ is one, while the entry $c_{kl}$ of $C$ shows the degree of closeness between nodes $k$ and $l$. We will see the influence of selecting $\{c_{kl}\}$ on the performance of the resulting algorithms in later simulations.

Similarly, we introduce an $N\times N$ adjacency matrix $\Phi$ with the element $\varphi_{kl}\in\{0,1\}$, in which $\varphi_{kl}=1$ if node $k$ is linked with node $l$; otherwise 0.

We also introduce the extended  aggregate matrix $G$
$$G\triangleq C\otimes I_{M},\;\;\;(NM\times NM)$$
where $\otimes$ is the Kronecker product operation and $I_{M}$ is the $M\times M$ identity matrix.

\subsection{Temporal-spatial recursion relation}

The temporal-spatial relation across network need to be considered as a starting point of convergence analysis. First, the diffusion strategy leads to the frequent spatial interaction between the neighborhoods, thereby each node $k$ is influenced by both local information such as $f_{k}$ and spatial information from neighbours $l\in \mathcal{N}_{k}$ such as $\{f_{l},x_{l}\}$. Second, the iteration way decides that the estimates and the local collected information on each node $k$ are time-variant, i.e., $\{f^{i}_{k},x^{i}_{k}\}$.

To begin with (\ref{eq9}), we have

\begin{equation}\label{eq15}
\mathcal{X}_{G}^{i}=Gx_{G}^{i}.
\end{equation}

Using (\ref{eq15}), we rewrite the local diffusion GN update step (\ref{eq11}) as a global representation

\begin{equation}\label{eq16}
x_{G}^{i+1}=Gx^{i}_{G}-\alpha D_{G}^{i}.
\end{equation}
Accordingly, we get the global non-cooperative GN update step
\begin{equation}\label{eq17}
x_{G}^{i+1}=x^{i}_{G}-\alpha d_{G}^{i}.
\end{equation}

Subtracting $\overline{x}^{*}$ on both sides of the equation (\ref{eq16}) and embedding the equation (\ref{eq17}), we get

\begin{equation}\label{eq18}
x_{G}^{i+1}-\overline{x}^{*}=(Gx^{i}_{G}-x^{i}_{G})+(x^{i}_{G}-\overline{x}^{*}-\alpha d_{G}^{i})+\alpha(d_{G}^{i}-D_{G}^{i}).
\end{equation}

Using the triangle inequality for vectors, we get the following recursion

\begin{equation}\label{eq19}
\begin{aligned}
\|x_{G}^{i+1}-\overline{x}^{*}\|\leq\|Gx^{i}_{G}-x^{i}_{G}\|+\|x^{i}_{G}-\overline{x}^{*}-\alpha d^{i}\|\\
+\alpha\|D_{G}^{i}-d_{G}^{i}\|.
\end{aligned}
\end{equation}

The inequality (\ref{eq19}) can be regarded as a temporal-spatial recursion relation, where the superscript $i$ and the subscript $G$ reflect the evolution of diffusion GN algorithm from temporal and spatial dimensions, respectively. And we establish the relation between diffusion GN and non-cooperative diffusion algorithms from the global perspective.

For the first term of the right side of (\ref{eq19}), we have

\begin{equation}\label{eq20}
\begin{aligned}
\|Gx^{i}_{G}-x^{i}_{G}\|&=\|Gx^{i}_{G}-G\overline{x}^{*}+(\overline{x}^{*}-x^{i}_{G})\|\\
&\leq \|Gx^{i}_{G}-G\overline{x}^{*}\|+\|x^{i}_{G}-\overline{x}^{*}\|\\
&\leq \|G\|_{F}\|x^{i}_{G}-\overline{x}^{*}\|+\|x^{i}_{G}-\overline{x}^{*}\|\\
&= (\|G\|_{F}+1)\|x^{i}_{G}-\overline{x}^{*}\|,
\end{aligned}
\end{equation}
where we use $G\overline{x}^{*}=\overline{x}^{*}$ based on the property of $G$.

For the second term of the right side of (\ref{eq19}), we have the following conclusion.

\textbf{Lemma 1.} Let Assumptions 1 and 2 hold. The norm of global vector $x^{i}_{G}-\overline{x}^{*}-\alpha d^{i}$ satisfies the following recursion
\begin{equation}\label{eq89}
\begin{aligned}
&\|x^{i}_{G}-\overline{x}^{*}-\alpha d^{i}\|\leq t_{1} \|x^{i}_{G}-\overline{x}^{*}\|^{2}+t_{2}\|x^{i}_{G}-\overline{x}^{*}\|,
\end{aligned}
\end{equation}
where
\begin{equation}\label{eq899}
\begin{aligned}
t_{1}\triangleq\dfrac{\alpha\omega}{2\Sigma_{min}},t_{2}\triangleq\frac{(1-\alpha)\Sigma_{max}}{\Sigma_{min}}+\frac{\sqrt{2}N\alpha\omega e_{min}}{\Sigma_{min}^{2}}
\end{aligned}
\end{equation}

\textit{Proof}: See Appendix A.

Given (\ref{eq20}) and (\ref{eq89}), we rewrite the temporal-spatial recursion relation (\ref{eq19}) as
\begin{equation}\label{eq27}
\begin{aligned}
\|x_{G}^{i+1}-\overline{x}^{*}\|&\leq t_{1}\|x^{i}_{G}-\overline{x}^{*}\|^{2}\\
&+(t_{2}+\|G\|_{F}+1)\|x^{i}_{G}-\overline{x}^{*}\|+\alpha\|D_{G}^{i}-d_{G}^{i}\|.
\end{aligned}
\end{equation}

Given the above, the left side of (\ref{eq27}) is the network deviation at time $i+1$, while the right side of (\ref{eq27}) will be related to the network deviation $\|x^{i}_{G}-\overline{x}^{*}\|$ at time $i$ if we can confirm that $\|D_{G}^{i}-d_{G}^{i}\|$ shares the same character or is  bounded by a given constant $\xi$. Then, we can establish the relation of the network deviation between the successive two times in diffusion GN.

\subsection{Boundness of descent discrepancy}
$D_{G}^{i}-d_{G}^{i}$ denotes the GN descent discrepancy over network between two modes of cooperative and non-cooperative. To decide the boundness of the discrepancy, we first evaluate the entry of $D_{G}^{i}-d_{G}^{i}$, i.e., $D_{k}^{i}-d_{k}^{i}$.

To begin the process, we write the entry as
\begin{equation}\label{eq28}
\begin{aligned}
D_{k}^{i}-d_{k}^{i}=[Q_{k}^{i}(\mathcal{X})]^{-1}q_{k}^{i}(\mathcal{X})-[Q_{k}^{i}(x_{k}^{i})]^{-1}q_{k}^{i}(x_{k}^{i}),\;\;\;k\in\mathcal{N}.
\end{aligned}
\end{equation}

Because of the matrix inverse operator, we introduce two  quantities
\begin{equation}\label{eq288}
\begin{aligned}
S_{k}^{i}\triangleq Q_{k}^{i}(\mathcal{X})-Q_{k}^{i}(x_{k}^{i})\\
s_{k}^{i}\triangleq q_{k}^{i}(\mathcal{X})-q_{k}^{i}(x_{k}^{i}).
\end{aligned}
\end{equation}

And in order to lower the impact of inverse operator for our analysis, the known matrix expansion formula \cite{Johnson1991Topics} will be used frequently in our analysis. That is
\begin{equation}\label{eq299}
\begin{aligned}
(Z+\delta Z)^{-1}=\sum\limits_{u=0}^{\infty}(-1)^{u}(Z^{-1}\delta Z)^{u}Z^{-1}
\end{aligned}
\end{equation}
for any matrix $Z$ and $\delta Z$ if $\|Z^{-1}\delta Z\|<1$.

From (\ref{eq9}), $\mathcal{X}_{k}^{i}$ is a convex combination of $\{x_{l}^{i}\}$ for $l\in \mathcal{N}_{k}$. Thus, Assumptions 1 and 2 hold for $\mathcal{X}_{k}^{i}$.

Then we have
\begin{equation}\label{eq29}
\begin{aligned}
\|S_{k}^{i}\|&=\|\sum\limits_{l\in \mathcal{N}_{k}}[F^{T}_{l}(\mathcal{X}_{l}^{i})F_{l}(\mathcal{X}_{l}^{i})-F^{T}_{l}(x_{k}^{i})F_{l}(x_{k}^{i})]\|\\
&\leq \sum\limits_{l\in \mathcal{N}_{k}}\gamma_{F}\|\mathcal{X}_{l}^{i}-x_{k}^{i}\|
\end{aligned}
\end{equation}

and

\begin{equation}\label{eq30}
\begin{aligned}
\|s_{k}^{i}\|&=\|\sum\limits_{l\in \mathcal{N}_{k}}[F^{T}_{l}(\mathcal{X}_{l}^{i})f_{l}(\mathcal{X}_{l}^{i})-F^{T}_{l}(x_{k}^{i})f_{l}(x_{k}^{i})]\|\\
&\leq \sum\limits_{l\in \mathcal{N}_{k}}\gamma_{f}\|\mathcal{X}_{l}^{i}-x_{k}^{i}\|.
\end{aligned}
\end{equation}

From (\ref{eq29}) and (\ref{eq30}), both $\|S_{k}^{i}\|$ and $\|s_{k}^{i}\|$ depend on $\|\mathcal{X}_{l}^{i}-x_{k}^{i}\|$. We now study the boundness of $\mathcal{X}_{l}^{i}-x_{k}^{i}$. Before that, we define a $1\times N$ vector
$$c_{l}\triangleq row\{c_{l1},c_{l2},\ldots, c_{lN}\},\;\;l\in\mathcal{N}$$
which is the $l$ row of matrix $C$.

Evaluating the norm of $\mathcal{X}_{l}^{i}-x_{k}^{i}$, we get
\begin{equation}\label{eq31}
\begin{aligned}
\|\mathcal{X}_{l}^{i}-x_{k}^{i}\|&=\|c_{l}x_{G}^{i}-c_{l}\textmd{1}_{N}x_{k}^{i}\|\\
&\leq \|c_{l}\|\|x_{G}^{i}-\textmd{1}_{N}x_{k}^{i}\|\\
&\leq \|x_{G}^{i}-\textmd{1}_{N}x_{k}^{i}\|.
\end{aligned}
\end{equation}

The block quantity $x_{G}^{i}-\textmd{1}_{N}x_{k}^{i}$ represents the estimate difference across the network at time $i$ and is written by
$$x_{G}^{i}-\textmd{1}_{N}x_{k}^{i}=\textrm{col}\{x_{1}^{i}-x_{k}^{i},x_{2}^{i}-x_{k}^{i},\ldots,x_{N}^{i}-x_{k}^{i}\}$$
whose individual entry is a $M\times 1$ vector.

For the norms of $x_{l}^{i}-x_{k}^{i}$ and $\mathcal{X}_{l}^{i}-x_{k}^{i}$, $l,k\in \mathcal{N}$ and $i\geq1$, we have the following Lemmas.

\textbf{Lemma 2.} Let Assumptions 1 and 2 hold. The estimate difference between nodes $l$ and $k$ through the non-cooperative GN update (\ref{eq111}) is bounded by
\begin{equation}\label{eq32}
\begin{aligned}
\|x_{l}^{i}-x_{k}^{i}\|\leq \Pi^{i},\;\;i\geq1
\end{aligned}
\end{equation}
where
\begin{equation}\label{eq6333}
\begin{aligned}
\Pi^{i}\triangleq a_{2}\sum\limits_{j=1}^{i}(a_{1})^{j-1},
\end{aligned}
\end{equation}
\begin{equation}\label{eq58}
\begin{aligned}
a_{1}\triangleq 1+\frac{\alpha n_{kl}+2\alpha n_{kl}\gamma_{f}}{2n_{l}\sigma_{min}^{2}},
\end{aligned}
\end{equation}
\begin{equation}\label{eq59}
\begin{aligned}
a_{2}\triangleq \frac{ (n_{l}+3n_{k|l}+3n_{l|k})\alpha\sigma_{max}\varepsilon_{max}}{2n_{l}\sigma_{min}^{2}},
\end{aligned}
\end{equation}
$n_{kl}$ denotes the number of nodes that are both in $\mathcal{N}_{k}$ and $\mathcal{N}_{l}$, $n_{k|l}$ denotes the number of nodes that are in $\mathcal{N}_{k}$ and not in $\mathcal{N}_{l}$.

\textit{Proof}: See Appendix B.

\textbf{Lemma 3.} Let Assumptions 1 and 2 hold. The estimate difference between nodes $l$ and $k$ through the diffusion  GN update (\ref{eq11}) is bounded by
\begin{equation}\label{eq62}
\begin{aligned}
\|\mathcal{X}_{l}^{i}-x_{k}^{i}\|\leq N\Pi^{i},\;\;i\geq 1,
\end{aligned}
\end{equation}
and
\begin{equation}\label{eq69}
\begin{aligned}
\|[Q_{k}^{i}(x_{k}^{i})]^{-1}S_{k}^{i}\|<1
\end{aligned}
\end{equation}
always holds under the sufficient condition
\begin{equation}\label{eq70}
\begin{aligned}
n_{kl}>0,
\end{aligned}
\end{equation}
where $a_{1}$, $a_{2}$, $Q_{k}^{i}(x_{k}^{i})$ and $S_{k}^{i}$ are assigned by (\ref{eq58}), (\ref{eq59}), (\ref{eq121}) and (\ref{eq288}), respectively.

\textit{Proof}: See Appendix D.

The condition (\ref{eq70}) means that any two nodes $k$ and $l$ in the network have at least one common neighboring node, which is more likely to be achieved by a small and dense
network. However, the condition can be relaxed in practice by allowing that all nodes are linked over single-hop or multi-hops so that it holds for the large scale networks. Thus, it is reasonable that the sufficient condition  $\|[Q_{k}^{i}(x_{k}^{i})]^{-1}S_{k}^{i}\|<1$ for applying the  expansion formula in (\ref{eq28}) always holds under Lemma 2.

Thus, we use the expansion formula (\ref{eq299}) and the norm operator on (\ref{eq28}) as follows
\begin{equation}\label{eq71}
\begin{aligned}
&\|D_{k}^{i}-d_{k}^{i}\|\\
&=\|[Q_{k}^{i}(x_{k}^{i})+S_{k}^{i}]^{-1}[q_{k}^{i}(x_{k}^{i})+s_{k}^{i}]-[Q_{k}^{i}(x_{k}^{i})]^{-1}q_{k}^{i}(x_{k}^{i})\|\\
&=\|[\sum\limits_{u=0}^{\infty}(-1)^{u}((Q_{k}^{i}(x_{k}^{i}))^{-1}S_{k}^{i})^{u}][Q_{k}^{i}(x_{k}^{i})]^{-1}[q_{k}^{i}(x_{k}^{i})+s_{k}^{i}]\\
&\;\;\;\;\;-[Q_{k}^{i}(x_{k}^{i})]^{-1}q_{k}^{i}(x_{k}^{i})\|\\
&=\|[\sum\limits_{u=1}^{\infty}(-1)^{u}((Q_{k}^{i}(x_{k}^{i}))^{-1}S_{k}^{i})^{u}][Q_{k}^{i}(x_{k}^{i})]^{-1}[q_{k}^{i}(x_{k}^{i})+s_{k}^{i}]\\
&\;\;\;\;\;+[Q_{k}^{i}(x_{k}^{i})]^{-1}s_{k}^{i}\|\\
&\leq \|[Q_{k}^{i}(x_{k}^{i})]^{-1}s_{k}^{i}\|+\\
&\|[\sum\limits_{u=1}^{\infty} ((Q_{k}^{i}(x_{k}^{i}))^{-1}S_{k}^{i})^{u}]\| \|[Q_{k}^{i}(x_{k}^{i})]^{-1}\| (\|q_{k}^{i}(x_{k}^{i})\|+\|s_{k}^{i}\|)\\
&\leq \frac{N\gamma_{f}\Pi^{i}}{\sigma_{min}^{2}}+\frac{(\sigma_{max}\varepsilon_{max}+N\gamma_{f}\Pi^{i})\zeta_{i}}{\sigma_{min}^{2}(1-\zeta_{i})},
\end{aligned}
\end{equation}
where the last equality comes from the obtained results including (\ref{eq40}) (\ref{eq411}) (\ref{eq64}) (\ref{eq65}) and the definitions (\ref{eq6333}) (\ref{eq666}) of $\Pi^{i}$ and $\zeta_{i}\in (0,1)$ (see Appendixes C and D). From (\ref{eq522}), we know that $\Pi^{i}$ is a bounded quantity that depends on the network topology.

Finally, we obtain the boundness conclusion as follows
\begin{equation}\label{eq72}
\begin{aligned}
&\|D_{G}^{i}-d_{G}^{i}\|\leq N\|D_{k}^{i}-d_{k}^{i}\|\\
&\leq \frac{N^{2}\gamma_{f}\Pi^{i}}{\sigma_{min}^{2}}+\frac{(N\sigma_{max}\varepsilon_{max}+N^{2}\gamma_{f}\Pi^{i})\zeta_{i}}{\sigma_{min}^{2}(1-\zeta_{i})}\triangleq \xi.
\end{aligned}
\end{equation}

\subsection{Convergence with sufficient conditions }

Given the constant $\xi>0$ that satisfies  (\ref{eq72}), we rewrite the global recursion relation (\ref{eq27}) as
\begin{equation}\label{eq73}
\begin{aligned}
&\|x_{G}^{i+1}-\overline{x}^{*}\|\\
&\leq t_{1}\|G\|^{2}\|x^{i}_{G}-\overline{x}^{*}\|^{2}+t_{2}\|G\|\|x^{i}_{G}-\overline{x}^{*}\|+\alpha\xi,
\end{aligned}
\end{equation}
which can be regarded as a nonlinear discrete dynamical system. Let $y^{i}\triangleq \|x^{i}_{G}-\overline{x}^{*}\|$, we will simplify notation of (\ref{eq73}) with the general form
\begin{equation}\label{eq74}
\begin{aligned}
y^{i+1} \leq t_{1}\|G\|^{2}(y^{i})^{2}+ t_{2} \|G\| y^{i}+\alpha\xi,
\end{aligned}
\end{equation}
whose steady-state equilibrium is a level \cite{Galor2007Discrete} that solves
\begin{equation}\label{eq75}
\begin{aligned}
y =\phi(y)= t_{1}\|G\|^{2}y^{2}+ t_{2} \|G\| y+\alpha\xi.
\end{aligned}
\end{equation}
Note that the steady-state equilibrium means that the variable $y^{i}$ is invariant under the law of motion indicated by the dynamical system. With the expression (\ref{eq75}), it is easy to know that the recursion (\ref{eq74}) is governed by the dynamical system $y^{i+1} =\phi(y^{i})$. Thus, guaranteeing the stability of system $y^{i+1} =\phi(y^{i})$ will be needed.

Solving (\ref{eq75}), we get two steady-state equilibrium points as follows
\begin{equation}\label{eq76}
\begin{aligned}
y_{max} =\frac{1-t_{2} \|G\|+\sqrt{(1-t_{2} \|G\|)^{2}-4t_{1}\alpha\xi\|G\|^{2}}}{2t_{1}\|G\|^{2}}
\end{aligned}
\end{equation}
and
\begin{equation}\label{eq77}
\begin{aligned}
y_{min} =\frac{1-t_{2} \|G\|-\sqrt{(1-t_{2} \|G\|)^{2}-4t_{1}\alpha\xi\|G\|^{2}}}{2t_{1}\|G\|^{2}}
\end{aligned}
\end{equation}
with the condition
\begin{equation}\label{eq78}
\begin{aligned}
(1-t_{2} \|G\|)^{2}-4t_{1}\alpha\xi\|G\|^{2}\geq 0.
\end{aligned}
\end{equation}

The equilibrium point  of the dynamical system (\ref{eq74}) is locally stable if and only if \cite{Galor2007Discrete}
\begin{equation}\label{eq79}
\begin{aligned}
\Big|\frac{d(\phi(y))}{dy}\Big|<1,
\end{aligned}
\end{equation}
where $\frac{d(\phi(y))}{dy}$ is the first order derivative of $\phi(y)$ with respect to $y$.

Thus, we know that $y_{max}$ is unstable since
\begin{equation}\label{eq80}
\begin{aligned}
\Big| \frac{d(\phi(y))}{dy}\big|_{y_{max}}\Big|=\Big|1+\sqrt{(1-t_{2} \|G\|)^{2}-4t_{1}\alpha\xi\|G\|^{2}}\Big|>1,
\end{aligned}
\end{equation}
while $y_{min}$ can be stable if
\begin{equation}\label{eq81}
\begin{aligned}
\Big| \frac{d(\phi(y))}{dy}\big|_{y_{min}}\Big|=\Big|1-\sqrt{(1-t_{2} \|G\|)^{2}-4t_{1}\alpha\xi\|G\|^{2}}\Big|<1
\end{aligned}
\end{equation}
holds.

From (\ref{eq81}), we get the following constraints
\begin{equation}\label{eq82}
\begin{aligned}
\frac{1}{t_{2}+2\sqrt{t_{1}\alpha\xi}}<\|G\|< \frac{t_{2}+2\sqrt{t_{2}^{2}-t_{1}\alpha\xi}}{t_{2}^{2}-4t_{1}\alpha\xi}
\end{aligned}
\end{equation}
and
\begin{equation}\label{eq91}
\begin{aligned}
\textrm{max}\{&\frac{t_{2}^{2}\|G\|^{2}-2t_{2} \|G\|-3}{4t_{1}\xi\|G\|^{2}},0\}<\alpha\\
&<\textrm{min}\{\frac{t_{2}^{2}\|G\|^{2}-2t_{2} \|G\|+1}{4t_{1}\xi\|G\|^{2}},1\},
\end{aligned}
\end{equation}
where $t_{1}$ and $t_{2}$ are given by (\ref{eq899}).

According to the locally stable theory for a steady-state equilibria in discrete dynamical systems\cite{Galor2007Discrete}, under the conditions (\ref{eq82}) (\ref{eq91}), as long as the initial condition $y^{0}$ is smaller than $y_{max}$, the nonlinear system (\ref{eq74}) converges to the unique steady-state equilibrium point $y_{min}$. That is

\begin{equation}\label{eq911}
\begin{aligned}
\lim\limits_{i\rightarrow\infty}y^{i}=y_{min}, \;\;\textrm{if}\;y_{0}<y_{max}.
\end{aligned}
\end{equation}

Considering $\|x_{G}^{i}-\overline{x}^{*}\|\geq0$ and $y_{min}<0$,  we have
\begin{equation}\label{eq912}
\begin{aligned}
\lim\limits_{i\rightarrow\infty}\|x_{G}^{i}-\overline{x}^{*}\|=0, \;\;\textrm{if}\;\|x_{G}^{0}-\overline{x}^{*}\|<y_{max}.
\end{aligned}
\end{equation}
In other words, the ATU algorithm converges asymptotically to the minimizer $x^{*}$ if the initial global error $\|x_{G}^{0}-\overline{x}^{*}\|<y_{max}$ holds. Conversely, the initial condition $\|x_{G}^{0}-\overline{x}^{*}\|>y_{max}$ will lead to the instability of algorithm and the growing global error level.

\subsection{Convergence behaviors}

In this section, we try to provide a qualitative analysis of the convergence behaviors for ATU and non-cooperative algorithms. Starting from the global ATU update (\ref{eq16}) and non-cooperative GN update (\ref{eq17}), and subtracting $\overline{x}^{*}$ on both sides of  (\ref{eq16}) and (\ref{eq17}), we get
\begin{equation}\label{eq94}
\begin{aligned}
x_{G}^{i+1}-\overline{x}^{*}&=Gx^{i}_{G}-\overline{x}^{*}-\alpha D_{G}^{i}\\
&=Gx^{i}_{G}-G\overline{x}^{*}-\alpha D_{G}^{i}+\alpha D_{G}^{*,i},
\end{aligned}
\end{equation}
and
\begin{equation}\label{eq944}
\begin{aligned}
x_{G}^{i+1}-\overline{x}^{*}&=x^{i}_{G}-\overline{x}^{*}-\alpha d_{G}^{i}\\
&=x^{i}_{G}-\overline{x}^{*}-\alpha d_{G}^{i}+\alpha d_{G}^{*,i},
\end{aligned}
\end{equation}
respectively, where we denote $D_{G}^{*,i}\triangleq \textrm{col}\{D_{1}^{*,i},D_{2}^{*,i},\ldots,D_{N}^{*,i}\}$, $d_{G}^{*,i}\triangleq \textrm{col}\{d_{1}^{*,i},d_{2}^{*,i},\ldots,d_{N}^{*,i}\}$ and $D_{k}^{*,i}\triangleq [Q_{k}^{i}(\mathcal{X})]^{-1}q_{k}^{i}(x^{*})=0$, $d_{k}^{*,i}\triangleq [Q_{k}^{i}(x_{k}^{i})]^{-1}q_{k}^{i}(x^{*})=0$, since $q_{k}^{i}(x^{*})=0$ and $G\overline{x}^{*}=\overline{x}^{*}$.

Applying the Triangle Inequality on the norm of (\ref{eq94}) and (\ref{eq944}), we have
\begin{equation}\label{eq95}
\begin{aligned}
\|x_{G}^{i+1}-\overline{x}^{*}\|&\leq \|Gx^{i}_{G}-G\overline{x}^{*}\|+\alpha \|D_{G}^{i}- D_{G}^{*,i}\|\\
&=\|Gx^{i}_{G}-G\overline{x}^{*}\|+\alpha \|\Lambda_{D} \rho_{D}\|\\
&=\|Gx^{i}_{G}-G\overline{x}^{*}\|+\alpha \|\Lambda_{D}\|\| \rho_{D}\|,\\
\end{aligned}
\end{equation}
and \begin{equation}\label{eq955}
\begin{aligned}
\|x_{G}^{i+1}-\overline{x}^{*}\|&\leq \|x^{i}_{G}-\overline{x}^{*}\|+\alpha \|d_{G}^{i}- d_{G}^{*,i}\|\\
&=\|x^{i}_{G}-\overline{x}^{*}\|+\alpha \|\Lambda_{d} \rho_{d}\|\\
&=\|x^{i}_{G}-\overline{x}^{*}\|+\alpha \|\Lambda_{d}\|\| \rho_{d}\|,\\
\end{aligned}
\end{equation}
where we introduce the matrices
\begin{equation}\label{eq96}
\begin{aligned}
\Lambda_{D}\triangleq\textrm{diag}\{[Q_{1}^{i}(\mathcal{X})]^{-1},\ldots,[Q_{N}^{i}(\mathcal{X})]^{-1}\},\;(NM\times NM),
\end{aligned}
\end{equation}
\begin{equation}\label{eq966}
\begin{aligned}
\Lambda_{d}\triangleq\textrm{diag}\{[Q_{1}^{i}(x_{1}^{i})]^{-1},\ldots,[Q_{N}^{i}(x_{N}^{i})]^{-1}\},\;(NM\times NM),
\end{aligned}
\end{equation}
and the vectors
\begin{equation}\label{eq97}
\begin{aligned}
\rho_{D}\triangleq\textrm{col}\{q_{1}^{i}(\mathcal{X})-q_{1}^{i}(x^{*}),\ldots,q_{N}^{i}(\mathcal{X})-q_{N}^{i}(x^{*})\},\;(NM\times1),
\end{aligned}
\end{equation}
\begin{equation}\label{eq977}
\begin{aligned}
\rho_{d}\triangleq\textrm{col}\{q_{1}^{i}(x_{1}^{i})-q_{1}^{i}(x^{*}),\ldots,q_{N}^{i}(x_{N}^{i})-q_{N}^{i}(x^{*})\},\;(NM\times1).
\end{aligned}
\end{equation}
Then, $\| \rho_{D}\|$ and $\| \rho_{d}\|$ are bounded as
\begin{equation}\label{eq98}
\begin{aligned}
\| \rho_{D}\|&=[\sum\limits_{k=1}^{N}\|\sum\limits_{l\in\mathcal{N}_{k}}(F^{T}_{l}(\mathcal{X}_{l}^{i})f_{l}(\mathcal{X}_{l}^{i})-F^{T}_{l}(x^{*})f_{l}(x^{*})\|^{2}]^{\frac{1}{2}}\\
&\leq [\sum\limits_{k=1}^{N}\sum\limits_{l\in\mathcal{N}_{k}}\|(F^{T}_{l}(\mathcal{X}_{l}^{i})f_{l}(\mathcal{X}_{l}^{i})-F^{T}_{l}(x^{*})f_{l}(x^{*})\|^{2}]^{\frac{1}{2}}\\
&\leq [\sum\limits_{k=1}^{N}\sum\limits_{l\in\mathcal{N}_{k}}\gamma_{f}^{2}\|\mathcal{X}_{l}^{i}-x^{*}\|^{2}]^{\frac{1}{2}}\\
&=\gamma_{f}\|\Omega(\mathcal{X}_{G}^{i}-\overline{x}^{*})\|\\
\end{aligned}
\end{equation}
and\begin{equation}\label{eq988}
\begin{aligned}
\| \rho_{d}\|&=[\sum\limits_{k=1}^{N}\|\sum\limits_{l\in\mathcal{N}_{k}}(F^{T}_{l}(x_{l}^{i})f_{l}(x_{l}^{i})-F^{T}_{l}(x^{*})f_{l}(x^{*})\|^{2}]^{\frac{1}{2}}\\
&\leq [\sum\limits_{k=1}^{N}\sum\limits_{l\in\mathcal{N}_{k}}\|(F^{T}_{l}(x_{l}^{i})f_{l}(x_{l}^{i})-F^{T}_{l}(x^{*})f_{l}(x^{*})\|^{2}]^{\frac{1}{2}}\\
&\leq [\sum\limits_{k=1}^{N}\sum\limits_{l\in\mathcal{N}_{k}}\gamma_{f}^{2}\|x_{l}^{i}-x^{*}\|^{2}]^{\frac{1}{2}}\\
&=\gamma_{f}\|\Omega(x_{G}^{i}-\overline{x}^{*})\|,\\
\end{aligned}
\end{equation}
where $\Omega$ is a $NN\times N$ matrix that can be written as a $N\times1$ block vector whose $k$ entry is the diagonal matrix $\textrm{diag}\{\varphi_{k1},\ldots,\varphi_{kN}\}$.

Substituting (\ref{eq98}) and (\ref{eq988}) into (\ref{eq95}) and (\ref{eq955}), respectively, we obtain the error recursions for our ATU algorithm and non-cooperative algorithm as follows
\begin{equation}\label{eq99}
\begin{aligned}
\|x_{G}^{i+1}-&\overline{x}^{*}\|\leq \|Gx^{i}_{G}-G\overline{x}^{*}\|+\alpha \gamma_{f}\|\Lambda_{D}\|\|\Omega(\mathcal{X}_{G}^{i}-\overline{x}^{*})\|\\
&\leq \|Gx^{i}_{G}-G\overline{x}^{*}\|+\alpha \gamma_{f}\|\Lambda_{D}\|\|\Omega\|\|Gx^{i}_{G}-G\overline{x}^{*}\|\\
&\leq (1+\alpha \gamma_{f}\|\Lambda_{D}\|\|\Omega\|)\|G (x^{i}_{G}-\overline{x}^{*})\|
\end{aligned}
\end{equation}
and
\begin{equation}\label{eq999}
\begin{aligned}
\|x_{G}^{i+1}-&\overline{x}^{*}\|\leq \|x^{i}_{G}-\overline{x}^{*}\|+\alpha \gamma_{f}\|\Lambda_{d}\|\|\Omega(x_{G}^{i}-\overline{x}^{*})\|\\
&\leq (1+\alpha \gamma_{f}\|\Lambda_{d}\|\|\Omega\|)\|x^{i}_{G}-\overline{x}^{*}\|.
\end{aligned}
\end{equation}

From (\ref{eq6}), we know that $\mathcal{X}_{k}^{i}$ is a convex combination of $\{x_{l}^{i}\}$ for $l\in \mathcal{N}_{k}$. Thus, Assumption 1  holds for $\mathcal{X}_{k}^{i}$. Under Assumption 1(1), we have
\begin{equation}\label{eq244}
\begin{aligned}
\|[Q_{k}^{i}(\mathcal{X})]^{-1}\|&=\|[\sum\limits_{l\in\mathcal{N}_{k}}F^{T}_{l}(\mathcal{X}_{l}^{i})F_{l}(\mathcal{X}_{l}^{i})]^{-1}\|\leq\frac{1}{n_{k}\sigma_{min}}
\end{aligned}
\end{equation}
and
\begin{equation}\label{eq255}
\|[Q_{k}^{i}(x^{i}_{k})]^{-1}\|\leq\frac{1}{n_{k}\sigma_{min}}.
\end{equation}
Furthermore, we have
\begin{equation}\label{eq266}
\begin{aligned}
\|\Lambda_{D}\|\leq\frac{1}{n_{k}\sigma_{min}}
\end{aligned}
\end{equation}
and
\begin{equation}\label{eq277}
\begin{aligned}
\|\Lambda_{d}\|\leq\frac{1}{n_{k}\sigma_{min}}.
\end{aligned}
\end{equation}
Thus, we know that $1+\alpha \gamma_{f}\|\Lambda_{D}\|\|\Omega\|$ and $1+\alpha \gamma_{f}\|\Lambda_{d}\|\|\Omega\|$ are
upper bounded by a small common constant when the small step size is selected.

The recursions (\ref{eq99}) and (\ref{eq999}) describe how the global error evolves over time for diffusion and non-cooperative GN algorithms, respectively.  It is important to note the difference between the linear structure (\ref{eq99}) and the nonlinear structure (\ref{eq73}). If we replace the lesser-or-equal  with an equal sign in (\ref{eq99}), the resulted linear system will be unstable due to $(1+\alpha \gamma_{f}\|\Lambda_{D}\|\|\Omega\|)\|G\|>1$ if $\|G\|\geq1$ \cite{Galor2007Discrete}.  However,  under guaranteed convergence conditions for  diffusion ATU, (\ref{eq99}) and (\ref{eq999}) reveal  the qualitative behavior of global error reduction in ATU and non-cooperative GN algorithms, respectively.

To analyze the convergence behavior of diffusion GN algorithm, we introduce the spectral radius of a square matrix, which is defined as the largest absolute value among its eigenvalues and denoted by $\rho(\cdot)$. Because of $G= C\otimes I_{M}$, we have
\begin{equation}\label{eq277-1}
\begin{aligned}
\rho(G)&=|\lambda_{max}(G)|=|\lambda_{max}(C\otimes I_{M})|\\
&\leq|\lambda_{max}(C)||\lambda_{max}(I_{M})|=|\lambda_{max}(C)|=1,
\end{aligned}
\end{equation}
where $|\lambda_{max}(C)|=1$ is based on the known conclusion (Appendix C of \cite{sayed2014diffusion}) if $C$ satisfies (\ref{eq14}). Thus, (\ref{eq277-1}) indicates that all eigenvalues of $G$ are smaller than 1, i.e., $|\lambda(G)|\leq1$.

Whenever we select an aggregation matrix $C$ based on (\ref{eq14}) so that $|\lambda(C)|\leq1$, the spectral radius of $G$ representing cooperative diffusion case is generally smaller than the spectral radius of $I_{N}$ representing non-cooperation case. That is, the cooperative diffusion GN algorithm will enforce a reduction of error $x_{G}^{i}-\overline{x}^{*}$ over the noncooperative version at every iteration. In other words, the error norm $\|x^{i}_{G}-\overline{x}^{*}\|$ in cooperation strategy decays more rapidly than that in non-cooperation strategy. The above analysis confirms the role of diffusion step in GN algorithm for improvements on convergence rate.

\section{Conclusions}

In this paper, we analyze the convergence of a cooperative  diffusion GN paradigm for nonlinear least squares problems in a distributed networked system. By the presented theoretical results, we show that cooperation strategy can obtain a room for improvement in term of convergence and guarantee algorithm's convergence when the derived sufficient conditions are satisfied, i.e., the good initial guesses, reasonable step size values and network connectivity. In order to avoid data incest and double counting, future works will focus on optimal estimate fusion, by which the network can adaptively adjust the weights or select the good nodes participating in the estimation.

\begin{appendices}

\section{Proof of Lemma 1}

The vector $x^{i}_{G}-\overline{x}^{*}-\alpha d_{G}^{i}$ is the global representation of $x_{k}^{i}-x^{*}-\alpha d_{k}^{i}$, which is written by
\begin{equation}\label{eq83}
\begin{aligned}
&x_{k}^{i}-x^{*}-\alpha d_{k}^{i}=x_{k}^{i}-x^{*}\\
&\;\;\;\;-\alpha[F^{T}_{l\in\mathcal{N}_{k}}(x_{k}^{i})F_{l\in\mathcal{N}_{k}}(x_{k}^{i})]^{-1}F^{T}_{l\in\mathcal{N}_{k}}(x_{k}^{i})f_{l\in\mathcal{N}_{k}}(x_{k}^{i})\\
&=x_{k}^{i}-x^{*}-\alpha[F_{l\in\mathcal{N}_{k}}(x_{k}^{i})]^{+}f_{l\in\mathcal{N}_{k}}(x_{k}^{i})\\
&=x_{k}^{i}-x^{*}-\alpha[F_{l\in\mathcal{N}_{k}}(x_{k}^{i})]^{+}f_{l\in\mathcal{N}_{k}}(x_{k}^{i})\\
&\;\;\;\;+\alpha[F_{l\in\mathcal{N}_{k}}(x^{*})]^{+}f_{l\in\mathcal{N}_{k}}(x^{*}),
\end{aligned}
\end{equation}
where $F^{+}(\cdot)$ denotes the generalized inverse of matrix $F(\cdot)$ and we use $[F_{l\in\mathcal{N}_{k}}(x^{*})]^{+}f_{l\in\mathcal{N}_{k}}(x^{*})=0$ according to Assumption 1.

According to the Assumption 1 and 2, we have the following inferences
\begin{equation}\label{eq84}
\begin{aligned}
\|A(x^{i}_{G})\|\leq \Sigma_{max},\;\;\;\|[A(x^{i}_{G})]^{+}\|\leq \frac{1}{\Sigma_{min}}
\end{aligned}
\end{equation}
and
\begin{equation}\label{eq85}
\begin{aligned}
\|A(x^{i}_{G})-A(\overline{x}^{*})\|\leq \omega \|x^{i}_{G} -\overline{x}^{*}\|.
\end{aligned}
\end{equation}

Thus, we can write the global representation as
\begin{equation}\label{eq86}
\begin{aligned}
&x^{i}_{G}-\overline{x}^{*}-\alpha d_{G}^{i}=[A(x^{i}_{G})]^{+}A(x^{i}_{G})(x^{i}_{G}-\overline{x}^{*})\\
&\;\;\;-\alpha [A(x^{i}_{G})]^{+}b(x^{i}_{G})+\alpha [A(\overline{x}^{*})]^{+}b(\overline{x}^{*})\\
&=[A(x^{i}_{G})]^{+}[A(x^{i}_{G})(x^{i}_{G}-\overline{x}^{*})-\alpha b(x^{i}_{G})+\alpha b(\overline{x}^{*})]\\
&\;\;\;+\alpha[(A(\overline{x}^{*}))^{+}-(A(x^{i}_{G}))^{+}]b(\overline{x}^{*})
\end{aligned}
\end{equation}

Using the mean-value theorem and (\ref{eq84}) (\ref{eq85}), we have
\begin{equation}\label{eq87}
\begin{aligned}
&\|\alpha b(\overline{x}^{*})-\alpha b(x^{i}_{G})-A(x^{i}_{G})(\overline{x}^{*}-x^{i}_{G})\|\\
&=\|\alpha \int_{0}^{1}A(x^{i}_{G}+u(\overline{x}^{*}-x^{i}_{G}))(\overline{x}^{*}-x^{i}_{G})du\\
&\;\;\;-A(x^{i}_{G})(\overline{x}^{*}-x^{i}_{G})\|\\
&=\|\alpha \int_{0}^{1}[A(x^{i}_{G}+u(\overline{x}^{*}-x^{i}_{G}))-A(x^{i}_{G})](\overline{x}^{*}-x^{i}_{G})du\\
&\;\;\;-(1-\alpha)A(x^{i}_{G})(\overline{x}^{*}-x^{i}_{G})\|\\
&\leq \alpha \int_{0}^{1}\|A(x^{i}_{G}+u(\overline{x}^{*}-x^{i}_{G}))-A(x^{i}_{G})\|du\|x^{i}_{G}-\overline{x}^{*}\|\\
&\;\;\;+(1-\alpha)\Sigma_{max}\|x^{i}_{G}-\overline{x}^{*}\|\\
&\leq \frac{\alpha\omega}{2} \|x^{i}_{G}-\overline{x}^{*}\|^{2}+(1-\alpha)\Sigma_{max}\|x^{i}_{G}-\overline{x}^{*}\|.\\
\end{aligned}
\end{equation}

Applying the Lemma 1 in \cite{Salzo2012Convergence}, we have
\begin{equation}\label{eq88}
\begin{aligned}
&\|(A(x^{i}_{G}))^{+}-(A(\overline{x}^{*}))^{+}\|\leq \\
&\sqrt{2} \|(A(\overline{x}^{*}))^{+}\|\|(A(x^{i}_{G}))^{+}\|\|A(x^{i}_{G})-A(\overline{x}^{*})\|\\
&\leq \frac{\sqrt{2}\omega}{\Sigma_{min}^{2}}\|x^{i}_{G}-\overline{x}^{*}\|.
\end{aligned}
\end{equation}
and
\begin{equation}\label{eq90}
\begin{aligned}
&\|b(\overline{x}^{*})\|\leq \sum\limits_{k=1}^{N}\|f_{l\in\mathcal{N}_{k}}(x^{*})\|=Ne_{min}.
\end{aligned}
\end{equation}

Therefore, substituting (\ref{eq87}) (\ref{eq88}) and (\ref{eq90}) into (\ref{eq86}) leads to (\ref{eq89}).

\section{Proof of Lemma 2}
To use the expansion formula (\ref{eq299}), we define
\begin{equation}\label{eq33}
\begin{aligned}
E_{kl}^{i}\triangleq Q_{k}^{i}(x_{k}^{i})-Q_{l}^{i}(x_{l}^{i})
\end{aligned}
\end{equation}

and
\begin{equation}\label{eq34}
\begin{aligned}
e_{kl}^{i}\triangleq q_{k}^{i}(x_{k}^{i})-q_{l}^{i}(x_{l}^{i}).
\end{aligned}
\end{equation}

Using (\ref{eq15}), we can rewrite (\ref{eq33}) and (\ref{eq34}) as
\begin{equation}\label{eq35}
\begin{aligned}
E_{kl}^{i}=\sum\limits_{u\in\mathcal{N}_{k}}F^{T}_{u}(x_{k}^{i})F_{u}(x_{k}^{i})-\sum\limits_{t\in\mathcal{N}_{l}}F^{T}_{t}(x_{l}^{i})F_{t}(x_{l}^{i})
\end{aligned}
\end{equation}
and
\begin{equation}\label{eq36}
\begin{aligned}
e_{kl}^{i}=\sum\limits_{u\in\mathcal{N}_{k}}F^{T}_{u}(x_{k}^{i})f_{u}(x_{k}^{i})-\sum\limits_{t\in\mathcal{N}_{l}}F^{T}_{t}(x_{l}^{i})f_{t}(x_{l}^{i}).
\end{aligned}
\end{equation}

Now we use the mathematical induction to obtain the results of Lemma 2.
\subsection{Initial Case: $i = 1$}
Given by the same initial estimate $x^{0}=x_{k}^{0}=x_{l}^{0}, k,l\in \mathcal{N}$ for all nodes in the network and (\ref{eq111}), we have

\begin{equation}\label{eq37}
\begin{aligned}
&x_{l}^{1}-x_{k}^{1}=\alpha[Q_{k}^{0}(x_{k}^{0})]^{-1}q_{k}^{0}(x_{k}^{0})-\alpha[Q_{l}^{0}(x_{l}^{0})]^{-1}q_{l}^{0}(x_{l}^{0})\\
&=\alpha[Q_{l}^{0}(x_{l}^{0})+E_{kl}^{0}]^{-1}(q_{l}^{0}(x_{l}^{0})+e_{kl}^{0})-\alpha[Q_{l}^{0}(x_{l}^{0})]^{-1}q_{l}^{0}(x_{l}^{0}).
\end{aligned}
\end{equation}

We now consider the conditions (the below Corollary 1) for applying the expansion formula (\ref{eq299}) when we let $Z=Q_{l}^{0}(x_{l}^{0})$ and $\delta Z=E_{kl}^{0}$.

\textbf{Corollary 1.} Let Assumptions 1 and 2 hold. The following recursion can be obtained
\begin{equation}\label{eq377}
\begin{aligned}
\|E_{kl}^{i}\|\leq n_{kl}\gamma_{F}\|x_{k}^{i}-x_{l}^{i}\|+( n_{k|l}+n_{l|k})\sigma_{max}^{2},
\end{aligned}
\end{equation}
and $\|[Q_{l}^{0}(x_{l}^{0})]^{-1}E_{kl}^{0}\|<1$ holds when the following condition
\begin{equation}\label{eq444}
\begin{aligned}
 \frac{n_{k|l}+n_{l|k}}{n_{l}}<\frac{\sigma_{min}^{2}}{\sigma_{max}^{2}}<1
\end{aligned}
\end{equation}
is satisfied.

\textit{Proof}: See Appendix C.

From (\ref{eq444}),  $\|[Q_{l}^{0}(x_{l}^{0})]^{-1}E_{kl}^{0}\|<1$ holds for a reasonable large denominator $n_{l}$ and a reasonable small numerator $n_{k|l}+n_{l|k}$. In other words, a high connectivity for the network is helpful for diffusion GN algorithm.

Under Corollary 1, we use the the expansion formula (\ref{eq299}) as follows
\begin{equation}\label{eq44}
\begin{aligned}
&[Q_{l}^{0}(x_{l}^{0})+E_{kl}^{0}]^{-1}[q_{l}^{0}(x_{l}^{0})+e_{kl}^{0}]-[Q_{l}^{0}(x_{l}^{0})]^{-1}q_{l}^{0}(x_{l}^{0})\\
&=[\sum\limits_{u=0}^{\infty}(-1)^{u}((Q_{l}^{0}(x_{l}^{0}))^{-1}E_{kl}^{0})^{u}((Q_{l}^{0}(x_{l}^{0}))^{-1}][q_{l}^{0}(x_{l}^{0})+e_{kl}^{0}]\\
&\;\;\;\;\;-[Q_{l}^{0}(x_{l}^{0})]^{-1}q_{l}^{0}(x_{l}^{0})\\
&=[\sum\limits_{u=1}^{\infty}(-1)^{u}((Q_{l}^{0}(x_{l}^{0}))^{-1}E_{kl}^{0})^{u}][(Q_{l}^{0}(x_{l}^{0}))^{-1}][q_{l}^{0}(x_{l}^{0})+e_{kl}^{0}]\\
&\;\;\;\;\;+[Q_{l}^{0}(x_{l}^{0})]^{-1}e_{kl}^{0}
\end{aligned}
\end{equation}

Substituting (\ref{eq44}) into (\ref{eq37}), for convenience of notation, we define
\begin{equation}\label{eq45}
\begin{aligned}
x_{l}^{1}-x_{k}^{1}\triangleq\alpha (p_{1}^{0}+p_{2}^{0})
\end{aligned}
\end{equation}
where $p_{1}^{0}$ and $p_{2}^{0}$ correspond sequentially to the last two terms of (\ref{eq44}), respectively.

Now we evaluate the norm of the vectors $p_{1}^{0}$ and $p_{2}^{0}$. According to the CBS inequality and the triangle inequality, we have
\begin{equation}\label{eq46}
\begin{aligned}
\|p_{1}^{0}\|&\leq (\sum\limits_{u=1}^{\infty}\|(Q_{l}^{0}(x_{l}^{0}))^{-1}E_{kl}^{0}\|)^{u}\|[Q_{l}^{0}(x_{l}^{0})]^{-1}\|\|q_{l}^{0}(x_{l}^{0})+e_{kl}^{0}\|\\
&\leq \frac{(n_{l}+n_{k|l}+n_{l|k})\sigma_{max}\varepsilon_{max}}{n_{l}\sigma_{min}^{2}}\sum\limits_{u=1}^{\infty}[\frac{(n_{k|l}+n_{l|k})\sigma_{max}^{2}}{n_{l}\sigma_{min}^{2}}]^{u}
\end{aligned}
\end{equation}
where we use (\ref{eq40}) (\ref{eq399}) (\ref{eq411}) for $i=0$ and (\ref{eq42}).

We set
$$\frac{(n_{k|l}+n_{l|k})\sigma_{max}^{2}}{n_{l}\sigma_{min}^{2}}=\mu_{kl}^{0}\in(0,1),$$ (\ref{eq46}) can be rewritten as
\begin{equation}\label{eq47}
\begin{aligned}
\|p_{1}^{0}\|&\leq \frac{(n_{l}+n_{k|l}+n_{l|k})\sigma_{max}\varepsilon_{max}\mu_{kl}^{0}}{n_{l}\sigma_{min}^{2}(1+\mu_{kl}^{0})}\\
&\leq \frac{(n_{l}+n_{k|l}+n_{l|k})\sigma_{max}\varepsilon_{max}}{2n_{l}\sigma_{min}^{2}}.
\end{aligned}
\end{equation}

For the norm of $p_{2}^{0}$, we have
\begin{equation}\label{eq48}
\begin{aligned}
\|p_{2}^{0}\|\leq \frac{(n_{k|l}+n_{l|k})\sigma_{max}\varepsilon_{max}}{n_{l}\sigma_{min}^{2}}.
\end{aligned}
\end{equation}

Therefore,
\begin{equation}\label{eq49}
\begin{aligned}
\|x_{l}^{1}-x_{k}^{1}\| \leq\alpha (\|p_{1}^{0}\|+\|p_{2}^{0}\|),
\end{aligned}
\end{equation}
where $\|p_{1}^{0}\|$ and $\|p_{2}^{0}\|$ are given by (\ref{eq48}) and (\ref{eq49}), respectively.

Therefore, we have $\|x_{l}^{1}-x_{k}^{1}\| \leq a_{2}$.

\subsection{Induction: $i=I$ and $i=I+1$}

For $i=I$ and any $l\neq k$,   let
\begin{equation}\label{eq50}
\begin{aligned}
\|x_{l}^{I}-x_{k}^{I}\|\leq a_{2}\sum\limits_{j=1}^{I}(a_{1})^{j-1}
\end{aligned}
\end{equation}
holds, where $a_{1}$ and $a_{2}$ are given by (\ref{eq58}) and (\ref{eq59}), respectively.

Then for $i=I+1$, we have
\begin{equation}\label{eq51}
\begin{aligned}
&x_{l}^{I+1}-x_{k}^{I+1}\\
&= x_{l}^{I}-x_{k}^{I}+\alpha[Q_{k}^{I}(x_{k}^{I})]^{-1}q_{k}^{I}(x_{k}^{I})-\alpha[Q_{l}^{I}(x_{l}^{I})]^{-1}q_{l}^{I}(x_{l}^{I})\\
&=x_{l}^{I}-x_{k}^{I}+\alpha[Q_{l}^{I}(x_{l}^{I})+E_{kl}^{I}]^{-1}(q_{l}^{I}(x_{l}^{I})+e_{kl}^{I})\\
&\;\;\;\;-\alpha[Q_{l}^{I}(x_{l}^{I})]^{-1}q_{l}^{I}(x_{l}^{I}).
\end{aligned}
\end{equation}

To apply the the expansion formula (\ref{eq299}) here, substituting (\ref{eq50}) into (\ref{eq41}) for $i=I$, the following condition
\begin{equation}\label{eq52}
\begin{aligned}
&\|[Q_{l}^{I}(x_{l}^{I})]^{-1}E_{kl}^{I}\|\\
&\leq \frac{n_{kl} \gamma_{F}a_{2}\sum\limits_{j=1}^{I}(a_{1})^{j-1}+(n_{k|l}+n_{l|k})\sigma_{max}^{2}}{n_{l}\sigma_{min}^{2}}\\
&\triangleq\mu_{kl}^{I}<1
\end{aligned}
\end{equation}
need to be satisfied. Substituting (\ref{eq58}) and (\ref{eq59}) into the above inequality, we get
\begin{equation}\label{eq53}
\begin{aligned}
&\alpha[(1+\alpha \theta)^{I+1}-(1+\alpha \theta)]\\
&<\frac{n_{l}\sigma_{min}^{2}-(n_{k|l}+n_{l|k})\sigma_{max}^{2}}{n_{kl} \gamma_{F}(n_{l}+3n_{k|l}+3n_{l|k})\sigma_{max}\varepsilon_{max}},
\end{aligned}
\end{equation}
where
\begin{equation}\label{eq533}
\begin{aligned}
\theta\triangleq \frac{n_{kl}+2n_{kl} \gamma_{f}}{2n_{l}\sigma_{min}^{2}}.
\end{aligned}
\end{equation}

The left side of (\ref{eq53}) is an increasing exponential function of the iteration time $I$. It is  a reasonable assumption that the inequality (\ref{eq53}) holds at any time $I$ when we set a sufficiently  small step size parameter $\alpha$. When (\ref{eq52}) is satisfied, we obtain the following useful conclusion
\begin{equation}\label{eq522}
\begin{aligned}
\|x_{l}^{I}-x_{k}^{I}\|\leq \frac{n_{l}\sigma_{min}^{2}-(n_{k|l}+n_{l|k})\sigma_{max}^{2}}{n_{kl} \gamma_{F}}.
\end{aligned}
\end{equation}

For (\ref{eq51}), we can bound $\|x_{l}^{I+1}-x_{k}^{I+1}\|$ by using the expansion formula (\ref{eq299})
\begin{equation}\label{eq54}
\begin{aligned}
\|x_{l}^{I+1}-x_{k}^{I+1}\|\leq \|x_{l}^{I}-x_{k}^{I}\|+\alpha (\|p_{1}^{I}\|+\|p_{2}^{I}\|),
\end{aligned}
\end{equation}
where we define $p_{1}^{I}$ and $p_{2}^{I}$ that are similar to (\ref{eq45}), and
\begin{equation}\label{eq55}
\begin{aligned}
&\|p_{1}^{I}\|\leq (\sum\limits_{u=1}^{\infty}\|(Q_{l}^{I}(x_{l}^{I}))^{-1}E_{kl}^{I}\|)^{u}\|[Q_{l}^{I}(x_{l}^{I})]^{-1}\|\|q_{l}^{I}(x_{l}^{I})+e_{kl}^{I}\|\\
&\leq \frac{(n_{l}+n_{k|l}+n_{l|k})\sigma_{max}\varepsilon_{max}+n_{kl}\|x_{l}^{I}-x_{k}^{I}\|}{n_{l}\sigma_{min}^{2}}\sum\limits_{u=1}^{\infty}[\mu_{kl}^{I}]^{u}\\
&\leq \frac{(n_{l}+n_{k|l}+n_{l|k})\sigma_{max}\varepsilon_{max}+n_{kl}\|x_{l}^{I}-x_{k}^{I}\|}{n_{l}\sigma_{min}^{2}}\frac{\mu_{kl}^{I}}{1+\mu_{kl}^{I}}\\
&\leq \frac{(n_{l}+n_{k|l}+n_{l|k})\sigma_{max}\varepsilon_{max}+n_{kl}\|x_{l}^{I}-x_{k}^{I}\|}{2n_{l}\sigma_{min}^{2}}
\end{aligned}
\end{equation}
and
\begin{equation}\label{eq56}
\begin{aligned}
&\|p_{2}^{I}\|\leq \frac{n_{kl} \gamma_{f}\|x_{l}^{I}-x_{k}^{I}\|+(n_{k|l}+n_{l|k})\sigma_{max}\varepsilon_{max}}{n_{l}\sigma_{min}^{2}}.
\end{aligned}
\end{equation}

Substituting (\ref{eq55}) (\ref{eq56}) into (\ref{eq54}) and rearranging them, we obtain
\begin{equation}\label{eq57}
\begin{aligned}
\|x_{l}^{I+1}-x_{k}^{I+1}\|\leq a_{1}\|x_{l}^{I}-x_{k}^{I}\|+a_{2}.
\end{aligned}
\end{equation}

Substituting (\ref{eq50}) into (\ref{eq57}), we can obtain
\begin{equation}\label{eq60}
\begin{aligned}
\|x_{l}^{I+1}-x_{k}^{I+1}\|&\leq a_{1}a_{2}\sum\limits_{j=1}^{I}(a_{1})^{j-1}+a_{2}\\
&=a_{2}(a_{1} \sum\limits_{j=1}^{I}(a_{1})^{j-1}+1)\\
&=a_{2}\sum\limits_{j=1}^{I+1}(a_{1})^{j-1}.
\end{aligned}
\end{equation}

\section{Proof of Corollary 1}

Starting from (\ref{eq35}), we get
\begin{equation}\label{eq38}
\begin{aligned}
E_{kl}^{i}&=\sum\limits_{u\in\mathcal{N}_{k}}F^{T}_{u}(x_{k}^{i})F_{u}(x_{k}^{i})-\sum\limits_{t\in\mathcal{N}_{l}}F^{T}_{t}(x_{l}^{i})F_{t}(x_{l}^{i})\\
&=\sum\limits_{u\in\mathcal{N}_{kl}}[F^{T}_{u}(x_{k}^{i})F_{u}(x_{k}^{i})-F^{T}_{u}(x_{l}^{i})F_{u}(x_{l}^{i})]\\
&+\sum\limits_{u\in\mathcal{N}_{k|l}}F^{T}_{u}(x_{k}^{i})F_{u}(x_{k}^{i})-\sum\limits_{t\in\mathcal{N}_{l|k}}F^{T}_{t}(x_{l}^{i})F_{t}(x_{l}^{i})\\
\end{aligned}
\end{equation}
and
\begin{equation}\label{eq388}
\begin{aligned}
e_{kl}^{i}&=\sum\limits_{u\in\mathcal{N}_{k}}F^{T}_{u}(x_{k}^{i})f_{u}(x_{k}^{i})-\sum\limits_{t\in\mathcal{N}_{l}}F^{T}_{t}(x_{l}^{i})f_{t}(x_{l}^{i})\\
&=\sum\limits_{u\in\mathcal{N}_{kl}}[F^{T}_{u}(x_{k}^{i})f_{u}(x_{k}^{i})-F^{T}_{u}(x_{l}^{i})f_{u}(x_{l}^{i})]\\
&+\sum\limits_{u\in\mathcal{N}_{k|l}}F^{T}_{u}(x_{k}^{i})f_{u}(x_{k}^{i})-\sum\limits_{t\in\mathcal{N}_{l|k}}F^{T}_{t}(x_{l}^{i})f_{t}(x_{l}^{i})
\end{aligned}
\end{equation}

Then from Assumption 2 and the CBS Inequality, we get
\begin{equation}\label{eq39}
\begin{aligned}
\|E_{kl}^{i}\|&\leq n_{kl} \gamma_{F}\|x_{k}^{i}-x_{l}^{i}\|+(n_{k|l}+n_{l|k})\sigma_{max}^{2}\\
\end{aligned}
\end{equation}
and
\begin{equation}\label{eq40}
\begin{aligned}
\|[Q_{l}^{i}(x_{l}^{i})]^{-1}\|\leq \frac{1}{n_{l}\sigma_{min}^{2}}.
\end{aligned}
\end{equation}

Based on Assumption 2, $\|f_{u\in\mathcal{N}_{k}}(x_{k}^{i})\|^{2}=\sum\limits_{u\in\mathcal{N}_{k}}\|f_{u}(x_{k}^{i})\|^{2}\leq e^{2}_{max}$ and $\|f_{u\in\mathcal{N}_{k}}(x_{k}^{i})\|^{2}\geq e_{min}^{2}$, thus $\|f_{u}(x_{k}^{i})\|$ has the upper and lower bounds for all $k\in\mathcal{N}$ and $x_{k}^{i}\in\mathbb{X}$. For convenience, we let $\varepsilon_{min}\leq\|f_{u}(x_{k}^{i})\|\leq\varepsilon_{max}$. Thus, we have $ \sigma_{min}\varepsilon_{min}\leq\|F_{u}^{T}(x_{k}^{i})f_{u}(x_{k}^{i})\|\leq\sigma_{max}\varepsilon_{max}$ and get
\begin{equation}\label{eq399}
\begin{aligned}
\|e_{kl}^{i}\| \leq n_{kl} \gamma_{f}\|x_{k}^{i}-x_{l}^{i}\|+(n_{k|l}+n_{l|k})\sigma_{max}\varepsilon_{max},
\end{aligned}
\end{equation}
and
\begin{equation}\label{eq411}
\begin{aligned}
\|q_{l}^{i}(x_{l}^{i})\|=\sum\limits_{t\in\mathcal{N}_{l}}F^{T}_{t}(x_{l}^{i})f_{t}(x_{l}^{i})\leq n_{l}\sigma_{max}\varepsilon_{max}.
\end{aligned}
\end{equation}

Thus, we have
\begin{equation}\label{eq41}
\begin{aligned}
\|[Q_{l}^{i}(x_{l}^{i})]^{-1}E_{kl}^{i}\|&\leq \frac{n_{kl} \gamma_{F}\|x_{k}^{i}-x_{l}^{i}\|+(n_{k|l}+n_{l|k})\sigma_{max}^{2}}{n_{l}\sigma_{min}^{2}}.
\end{aligned}
\end{equation}

To ensure $\|[Q_{l}^{0}(x_{l}^{0})]^{-1}E_{kl}^{0}\|<1$ when given $x_{k}^{0}=x_{l}^{0}$, we obtain the following condition
\begin{equation}\label{eq42}
\begin{aligned}
\|[Q_{l}^{0}(x_{l}^{0})]^{-1}E_{kl}^{0}\|\leq \frac{(n_{k|l}+n_{l|k})\sigma_{max}^{2}}{n_{l}\sigma_{min}^{2}}<1.
\end{aligned}
\end{equation}
To further rewrite  (\ref{eq42}), we get (\ref{eq444}).

\section{Proof of Lemma 3}

Under Lemma 2 and the inequality (\ref{eq31}), we know that $\|\mathcal{X}_{l}^{i}-x_{k}^{i}\|$ is also bounded by
\begin{equation}\label{eq633}
\begin{aligned}
\|\mathcal{X}_{l}^{i}-x_{k}^{i}\|\leq N\|x_{l}^{i}-x_{k}^{i}\|.
\end{aligned}
\end{equation}
Thus, (\ref{eq62}) can be obtained.

Thus, we can replace  (\ref{eq29}) (\ref{eq30}) as
\begin{equation}\label{eq63}
\begin{aligned}
\|S_{k}^{i}\|\leq Nn_{k}\gamma_{F}\Pi^{i}
\end{aligned}
\end{equation}
and
\begin{equation}\label{eq64}
\begin{aligned}
\|s_{k}^{i}\|\leq N n_{k}\gamma_{f}\Pi^{i}.
\end{aligned}
\end{equation}

Setting $Z=Q_{k}^{i}(x_{k}^{i})$ and $\delta Z=S_{k}^{i}$, we have
\begin{equation}\label{eq65}
\begin{aligned}
\|[Q_{k}^{i}(x_{k}^{i})]^{-1}S_{k}^{i}\|&\leq \|[Q_{k}^{i}(x_{k}^{i})]^{-1}\|\|S_{k}^{i}\|\\
&\leq \frac{Nn_{k}\gamma_{F}\Pi^{i}}{\sigma_{min}^{2}}.
\end{aligned}
\end{equation}

Now defining
\begin{equation}\label{eq666}
\begin{aligned}
\frac{Nn_{k}\gamma_{F}\Pi^{i}}{\sigma_{min}^{2}}\triangleq \zeta_{i}
\end{aligned}
\end{equation}
for subsequent use, and solving $\zeta_{i}<1$ and combining (\ref{eq522}), an inequality is obtained as
\begin{equation}\label{eq66}
\begin{aligned}
(Nn_{k}n_{l}-n_{kl})\sigma_{min}^{2}<Nn_{k}(n_{k|l}+n_{l|k})\sigma_{max}^{2}.
\end{aligned}
\end{equation}

By rewriting  (\ref{eq444}) as follows
\begin{equation}\label{eq67}
\begin{aligned}
(n_{k|l}+n_{l|k})\sigma_{max}^{2}<n_{l}\sigma_{min}^{2}
\end{aligned}
\end{equation}
and substituting it into  (\ref{eq66}) gives
\begin{equation}\label{eq68}
\begin{aligned}
n_{kl}>0.
\end{aligned}
\end{equation}

\end{appendices}

\bibliographystyle{IEEEtran}
\bibliography{mybibfile}

\begin{thebibliography}{10}
\providecommand{\url}[1]{#1}
\csname url@samestyle\endcsname
\providecommand{\newblock}{\relax}
\providecommand{\bibinfo}[2]{#2}
\providecommand{\BIBentrySTDinterwordspacing}{\spaceskip=0pt\relax}
\providecommand{\BIBentryALTinterwordstretchfactor}{4}
\providecommand{\BIBentryALTinterwordspacing}{\spaceskip=\fontdimen2\font plus
\BIBentryALTinterwordstretchfactor\fontdimen3\font minus
  \fontdimen4\font\relax}
\providecommand{\BIBforeignlanguage}[2]{{%
\expandafter\ifx\csname l@#1\endcsname\relax
\typeout{** WARNING: IEEEtran.bst: No hyphenation pattern has been}%
\typeout{** loaded for the language `#1'. Using the pattern for}%
\typeout{** the default language instead.}%
\else
\language=\csname l@#1\endcsname
\fi
#2}}
\providecommand{\BIBdecl}{\relax}
\BIBdecl

\bibitem{Schmidhuber2015Deep}
J.~Schmidhuber, ``Deep learning in neural networks: an overview.'' \emph{Neural
  Network}, vol.~61, pp. 85--117, 2015.

\bibitem{wang2018deep}
C.~C. Wang, K.~L. Tan, C.~T. Chen, Y.~H. Lin, S.~S. Keerthi, D.~Mahajan,
  S.~Sundararajan, and C.~J. Lin, ``Distributed newton methods for deep neural
  networks,'' \emph{Neural Computation}, vol.~30, no.~6, pp. 1673--1724, 2018.

\bibitem{Zhao2007Information}
T.~Zhao and A.~Nehorai, ``Information-driven distributed maximum likelihood
  estimation based on gauss-newton method in wireless sensor networks,''
  \emph{IEEE Transactions on Signal Processing}, vol.~55, no.~9, pp.
  4669--4682, 2007.

\bibitem{Xiao2013Convergence}
L.~Xiao and A.~Scaglione, ``Convergence and applications of a gossip-based
  gauss-newton algorithm,'' \emph{IEEE Transactions on Signal Processing},
  vol.~61, no.~21, pp. 5231--5246, 2013.

\bibitem{Kaur2017A}
A.~Kaur, P.~Kumar, and G.~P. Gupta, ``A novel dv-hop algorithm based on
  gauss-newton method,'' in \emph{Fourth International Conference on Parallel,
  Distributed and Grid Computing}, 2017, pp. 625--629.

\bibitem{Gratton2007Approximate}
S.~Gratton, A.~S. Lawless, and N.~K. Nichols, ``Approximate gauss¨cnewton
  methods for nonlinear least squares problems,'' \emph{Siam Journal on
  Optimization}, vol.~18, no.~1, pp. 106--132, 2007.

\bibitem{Bao2017Approximate}
J.~F. Bao, C.~Li, W.~P. Shen, J.~C. Yao, and S.~M. Guu, ``Approximate
  gauss-newton methods for solving underdetermined nonlinear least squares
  problems,'' \emph{Applied Numerical Mathematics}, vol. 111, pp. 92--110,
  2017.

\bibitem{Xiong2014Supervised}
X.~Xiong and F.~D.~L. Torre, ``Supervised descent method for solving nonlinear
  least squares problems in computer vision,'' \emph{Eprint Arxiv}, 2014.

\bibitem{Torre2008Parameterized}
F.~D.~L. Torre and M.~H. Nguyen, ``Parameterized kernel principal component
  analysis: Theory and applications to supervised and unsupervised image
  alignment,'' in \emph{Computer Vision and Pattern Recognition, IEEE
  Conference on}, 2008, pp. 1--8.

\bibitem{Schweiger2005Gauss}
M.~Schweiger, S.~R. Arridge, and I.~Nissila, ``Gauss-newton method for image
  reconstruction in diffuse optical tomography,'' \emph{Physics in Medicine \&
  Biology}, vol.~50, no.~10, pp. 2365--86, 2005.

\bibitem{souza2016target}
E.~L. Souza, E.~F. Nakamura, and R.~W. Pazzi, ``Target tracking for sensor
  networks: A survey,'' \emph{ACM Computing Surveys (CSUR)}, vol.~49, no.~2,
  p.~30, 2016.

\bibitem{bejar2011distributed}
B.~Bejar, P.~Belanovic, and S.~Zazo, ``Distributed consensus-based tracking in
  wireless sensor networks: a practical approach,'' in \emph{Signal Processing
  Conference, 19th European}.\hskip 1em plus 0.5em minus 0.4em\relax IEEE,
  2011, pp. 2019--2023.

\bibitem{Jensen2013Nonlinear}
J.~R. Jensen, M.~G. Christensen, and S.~H. Jensen, ``Nonlinear least squares
  methods for joint doa and pitch estimation,'' \emph{IEEE Transactions on
  Audio Speech \& Language Processing}, vol.~21, no.~5, pp. 923--933, 2013.

\bibitem{Zhuang2015distributed}
Y.~C.~J. Yong~Zhuang, Wei Sheng~Chin and C.~J. Lin, ``Distributed newton
  methods for regularized logistic regression,'' in \emph{2015 Pacific-Asia
  Conference on Knowledge Discovery and Data Mining}.\hskip 1em plus 0.5em
  minus 0.4em\relax Springer, 2015, pp. 690--703.

\bibitem{cosovic2016distributed}
M.~Cosovic and D.~Vukobratovic, ``Distributed gauss-newton method for ac state
  estimation: A belief propagation approach,'' in \emph{IEEE International
  Conference on Smart Grid Communications}, 2016.

\bibitem{Minot2016A}
A.~Minot, Y.~M. Lu, and N.~Li, ``A distributed gauss-newton method for power
  system state estimation,'' \emph{IEEE Transactions on Power Systems},
  vol.~31, no.~5, pp. 3804--3815, 2016.

\bibitem{Calafiore2010A}
G.~C. Calafiore, L.~Carlone, and M.~Wei, ``A distributed gauss-newton approach
  for range-based localization of multi agent formations,'' in \emph{IEEE
  International Symposium on Computer-Aided Control System Design}, 2010, pp.
  1152--1157.

\bibitem{Bejar2015Distributed}
B.~Bejar, P.~Belanovic, and S.~Zazo, ``Distributed consensus-based tracking in
  wireless sensor networks: A practical approach,'' in \emph{Signal Processing
  Conference, 2011 19th European}, 2011, pp. 2019--2023.

\bibitem{kelley1999iterative}
C.~T. Kelley, \emph{Iterative methods for optimization}.\hskip 1em plus 0.5em
  minus 0.4em\relax SIAM, 1999.

\bibitem{Bj1996Numerical}
A.~Bjorck, \emph{Numerical methods for least squares problems}.\hskip 1em plus
  0.5em minus 0.4em\relax SIAM, 1996.

\bibitem{Johnson1991Topics}
C.~R. Johnson, \emph{Topics in matrix analysis}.\hskip 1em plus 0.5em minus
  0.4em\relax Cambridge University Press,, 1991.

\bibitem{Eriksson2004Applied}
E.~D. Eriksson~Kenneth and J.~Claes, \emph{Applied Mathematics, Body and Soul:
  Derivatives and Geometry in \uppercase{IR}3}.\hskip 1em plus 0.5em minus
  0.4em\relax Springer, 2004.

\bibitem{Galor2007Discrete}
O.~Galor, \emph{Discrete Dynamical Systems}.\hskip 1em plus 0.5em minus
  0.4em\relax Springer, 2007.

\bibitem{sayed2014diffusion}
A.~H. Sayed, ``Diffusion adaptation over networks,'' in \emph{Academic Press
  Library in Signal Processing}.\hskip 1em plus 0.5em minus 0.4em\relax
  Elsevier, 2014, vol.~3, pp. 323--453.

\bibitem{Salzo2012Convergence}
S.~Salzo and S.~Villa, ``Convergence analysis of a proximal gauss-newton
  method,'' \emph{Computational Optimization \& Applications}, vol.~53, no.~2,
  pp. 557--589, 2012.

\end{thebibliography}

\end{document}